\documentclass{gtart}

\def\ifplaintex{\expandafter\ifx\csname documentclass\endcsname\relax}

\def\gtp{{\mathsurround=0pt\it $\cal G\mskip-2mu$eometry \&\ 
$\cal T\!\!$opology $\cal P\!$ublications}}  

\def\recd{{\small Received:\qua\receiveddate\ifx\reviseddate\relax
\else\qquad Revised:\qua\reviseddate\fi\par}} 

\def\lognumber#1{\def\thelognumber{#1}}
\def\volumenumber#1{\def\thevolumenumber{#1}}
\def\volumeyear#1{\def\thevolumeyear{#1}}
\def\papernumber#1{\def\thepapernumber{#1}}
\def\pagenumbers#1#2{\def\startpage{#1}\def\finishpage{#2}}
\def\published#1{\def\publishdate{#1}}

\def\received#1{\def\receiveddate{#1}}

\def\accepted#1{\def\accepteddate{#1}}
\def\asciititle#1{\def\theasciititle{#1}}
\def\covertitle#1{\def\thecovertitle{#1}}

\long\def\asciiabstract#1{\long\def\theasciiabstract{#1}}
\def\asciikeywords#1{\def\theasciikeywords{#1}}

\let\\\par\let\thelognumber\relax\let\thevolumenumber\relax
\let\thepapernumber\relax\let\thevolumeyear\relax\let\startpage\relax
\let\finishpage\relax\let\publishdate\relax\let\receiveddate\relax
\let\reviseddate\relax\let\accepteddate\relax\let\theasciititle\relax
\let\thecovertitle\relax\let\theasciiauthors\relax
\let\theasciiabstract\relax\let\theasciikeywords\relax

\let\theasciiemail\relax

\ifplaintex
\font\logobig=cmssbx10 scaled 3836
\font\logomed=cmssbx10 scaled 2557
\else
\font\logobig=cmssbx10 scaled 4200
\font\logomed=cmssbx10 scaled 2800
\fi

\long\def\makeagttitle{   
\count0=\startpage
\agt\hfill      
\hbox to 45truept{\vbox to 0pt{\vglue -13truept{\logomed A\kern -.37em{\logobig 
T}\kern -.38em G}\vss}\hss}
\break
{\small Volume \thevolumenumber\ (\thevolumeyear)
\startpage--\finishpage\nl
Published: \publishdate}

\vglue .25truein

{\parskip=0pt\leftskip 0pt plus
1fil\def\\{\par\smallskip}{\Large\bf\thetitle}\par\medskip} \vglue
0.05truein

{\parskip=0pt\leftskip 0pt plus 1fil\def\\{\par}{\sc\theauthors}
\par\medskip}%
 
\vglue 0.03truein 

{\small\leftskip 25truept\rightskip 25truept{\bf Abstract}\stdspace\theabstract

{\bf AMS Classification}\stdspace\theprimaryclass
\ifx\thesecondaryclass\relax\else; \thesecondaryclass\fi\par
{\bf Keywords}\stdspace \thekeywords\par}\vglue 7truept

}   

\ifplaintex
\font\phead=cmsl9 scaled 950
\font\pnum=cmbx10 scaled 913
\font\pfoot=cmsl9 scaled 950
\headline{\vbox to 0pt{\vskip -4.5mm\line{\small\phead\ifnum
\count0=\startpage ISSN 1472-2739 (on-line) 1472-2747 (printed)
\hfill {\pnum\folio}\else\ifodd\count0\def\\{ }%
\ifx\theshorttitle\relax\thetitle\else\theshorttitle\fi\hfill{\pnum\folio}
\else\def\\{ and }{\pnum\folio}\hfill\ifx\theshortauthors\relax\theauthors
\else\theshortauthors\fi\fi\fi}\vss}}
\footline{\vbox to 0pt{\vglue 0mm\line{\small\pfoot\ifnum\count0=\startpage
\copyright\ \gtp\hfill\else
\agt, Volume \thevolumenumber\ (\thevolumeyear)\hfill\fi}\vss}}
\else
\font=cmsl9 scaled 1050
\font\lnum=cmbx10 
\font=cmsl9 scaled 1050
\makeatletter
\def\@oddhead{{\small\lhead\ifnum\count0=\startpage ISSN 1472-2739 
(on-line) 1472-2747 (printed)\hfill {\lnum\number\count0}\else\ifodd\count0
\def\\{ }\ifx\theshorttitle\relax \thetitle \else\theshorttitle\fi\hfill
{\lnum\number\count0}\else\def\\{ and }{\lnum\number\count0}
\hfill\ifx\theshortauthors\relax 
\theauthors\else\theshortauthors\fi\fi\fi}}\def\@evenhead{\@oddhead}
\def\@oddfoot{\small\lfoot\ifnum\count0=\startpage\copyright\ \gtp\hfill\else
\agt, Volume \thevolumenumber\ (\thevolumeyear)\hfill\fi}
\def\@evenfoot{\@oddfoot}
\makeatother
\fi
\let\maketitlepage\makeagttitle
\let\makeshorttitle\maketitlepage
\let\maketitle\maketitlepage

\newwrite\gtoutfile
\long\gdef\makeheadfile{  
{\def\\{, }\def\s{ }
\immediate\openout\gtoutfile head.xxx
\immediate\write\gtoutfile{Proxy-for: \ifx\theasciiauthors\relax
\theauthors\else\theasciiauthors\fi\s<\ifx\theasciiemail\relax\theemail\else\theasciiemail\fi>}
\immediate\write\gtoutfile{\noexpand\\}
\immediate\write\gtoutfile{Authors: \ifx\theasciiauthors\relax
\theauthors\else\theasciiauthors\fi}
{\def\\{ }\immediate\write\gtoutfile{Title: \ifx\theasciititle\relax
\thetitle\else\theasciititle\fi}}
\immediate\write\gtoutfile{Subj-class: GT or SG, GR etc}
\immediate\write\gtoutfile{MSC-class: \theprimaryclass\ifx\thesecondaryclass\relax\else, \thesecondaryclass\fi}
\immediate\write\gtoutfile{Journal-ref: Algebr. Geom. Topol. \thevolumenumber\s
(\thevolumeyear) \startpage-\finishpage}
\immediate\write\gtoutfile{Comments: Published by Algebraic and
Geometric Topology at}
\immediate\write\gtoutfile{\s\s\s  http://www.maths.warwick.ac.uk/agt/AGTVol\thevolumenumber/agt-\thevolumenumber-\thepapernumber.abs.html}
\immediate\write\gtoutfile{\noexpand\\}
\immediate\write\gtoutfile{}
\ifx\theasciiabstract\relax
\immediate\write\gtoutfile{\theabstract}\else
\immediate\write\gtoutfile{\theasciiabstract}\fi
\immediate\write\gtoutfile{}
\immediate\write\gtoutfile{\noexpand\\}
\immediate\write\gtoutfile{}
\immediate\closeout\gtoutfile}}  

\def\maketitlepage{\makeagttitle\makeheadfile}
\let\makeshorttitle\maketitlepage
\let\maketitle\maketitlepage

\lognumber{40}
\volumenumber{3}
\volumeyear{2003}
\papernumber{40}
\published{10 November 2003}
\pagenumbers{1119}{1138}
\received{24 October 2003}
\accepted{5 November2003}

\usepackage{amsmath}
\usepackage{amssymb}

\theoremstyle{plain}
\newtheorem{theorem}{Theorem}[section]

\newtheorem{corollary}[theorem]{Corollary}

\newtheorem{lemma}[theorem]{Lemma}
\newtheorem{proposition}[theorem]{Proposition}

\theoremstyle{definition}

\newtheorem{definition}[theorem]{Definition}

\newtheorem{remark}[theorem]{Remark}
\theoremstyle{remark}
\numberwithin{equation}{section}

\nocolon

\begin{document}
\title[Structure of mod two symmetric algebra]{Global structure of the 
mod two symmetric algebra, 
\\ $H^{\ast}(BO;\mathbb{F}_{2})$, over the Steenrod Algebra}
\asciititle{Global structure of the 
mod two symmetric algebra, H^*(BO;F_2), over the Steenrod Algebra}
\covertitle{Global structure of the 
mod two symmetric algebra,\\$H^*(BO;{\noexpand\bf F}_{2})$, 
over the Steenrod Algebra}
\authors{David J. Pengelley\\Frank Williams}
\address{New Mexico State University\\
Las Cruces, NM 88003, USA}
\email{davidp@nmsu.edu, frank@nmsu.edu}

\begin{abstract}
The algebra $\mathcal{S}$ of symmetric invariants over the field with two
elements is an unstable algebra over the Steenrod algebra $\mathcal{A}$, and
is isomorphic to the mod two cohomology of $BO$, the classifying space for
vector bundles. We provide a minimal presentation for $\mathcal{S}$ in the
category of unstable $\mathcal{A}$-algebras, i.e., minimal generators and
minimal relations.

{}From this we produce minimal presentations for various unstable $\mathcal{A}%
$-algebras associated with the cohomology of related spaces, such as the
$BO(2^{m}-1)$ that classify finite dimensional vector bundles, and the
connected covers of $BO$. The presentations then show that certain of these
unstable $\mathcal{A}$-algebras coalesce to produce the Dickson algebras of
general linear group invariants, and we speculate about possible related
topological realizability.

Our methods also produce a related simple minimal $\mathcal{A}$-module
presentation of the cohomology of infinite dimensional real projective space,
with filtered quotients the unstable modules $\mathcal{F}\left(
2^{p}-1\right)  /\mathcal{A}\overline{\mathcal{A}}_{p-2}$, as described in an
independent appendix.
\end{abstract}

\asciiabstract{The algebra S of symmetric invariants over the field with two
elements is an unstable algebra over the Steenrod algebra A, and
is isomorphic to the mod two cohomology of BO, the classifying space for
vector bundles. We provide a minimal presentation for S in the
category of unstable A-algebras, i.e., minimal generators and
minimal relations.

From this we produce minimal presentations for various unstable
A-algebras associated with the cohomology of related spaces, such as
the BO(2^m-1) that classify finite dimensional vector bundles, and
the connected covers of BO. The presentations then show that certain
of these unstable A-algebras coalesce to produce the
Dickson algebras of general linear group invariants, and we speculate
about possible related topological realizability.

Our methods also produce a related simple minimal A-module
presentation of the cohomology of infinite dimensional real projective
space, with filtered quotients the unstable modules
F(2^p-1)/A bar{A}_{p-2}, as described in an independent appendix.}



\asciikeywords{Symmetric algebra, Steenrod algebra, unstable algebra,
classifying space, Dickson algebra, BO, real projective space.}

\maketitle

\section{Introduction}

We continue our study \cite{ppw} of invariant algebras as unstable algebras
over the Steenrod algebra $\mathcal{A}$ by proving a structure theorem for the
algebra $\mathcal{S}$ of symmetric invariants over the field $\mathbb{F}_{2}$.
The algebra $\mathcal{S}$ is isomorphic to the mod two cohomology of $BO$, the
classifying space for vector bundles \cite{milnor}, and we identify the two.
We also make several applications to the cohomology of related spaces, which
then reveal a relationship between $\mathcal{S}$ and the Dickson algebras
\cite{wilkerson}.

Our goal is to provide a minimal presentation for $\mathcal{S}=H^{\ast
}(BO;\mathbb{F}_{2})$ in the category of unstable $\mathcal{A}$-algebras
\cite{steenrodepstein}, beginning with a minimally presented generating
$\mathcal{A}$-module and then introducing a minimal set of $\mathcal{A}%
$-algebra relations. This reveals how a minimal set of $\mathcal{A}$-module
building blocks for $\mathcal{S}$ fit together in its $\mathcal{A}$-algebra
structure. In brief, our main result (Theorem \ref{T:mainthm}) is that
$\mathcal{S}=H^{\ast}(BO;\mathbb{F}_{2})$ is minimally presented in the
category of unstable $\mathcal{A}$-algebras as the free unstable $\mathcal{A}%
$-algebra on the two-power Stiefel-Whitney classes $w_{2^{k}}$ modulo
relations expressing the fact that, for each $i\leq k-2$, $Sq^{2^{i}}w_{2^{k}%
}$ differs from $Sq^{2^{k-1}}Sq^{2^{i}}w_{2^{k-1}}$ by a decomposable. (By
contrast, and at first seemingly paradoxically, we shall also see (Theorem
\ref{T:injects}) that while $\mathcal{S}$ is generated as an $\mathcal{A}%
$-algebra by $\left\{  w_{2^{k}}\co k\geq0\right\}  $, with relations linking the
resulting algebra generators, in fact the $\mathcal{A}$-submodule of
$\mathcal{S}$ generated by $\{w_{m}\co m\geq0\}$ is a free unstable $\mathcal{A}%
$-module on all the Stiefel-Whitney classes.)

We apply this structure theorem to characterize similarly the cohomology
images $B^{\ast}(n)$ for the connected covers of $BO$ (Theorem
\ref{T:conncover}) \cite{kochman}, which include the full cohomology algebras
of $BSO$, $BSpin$, and $BO\left\langle 8\right\rangle $. We likewise
characterize the quotients $H^{\ast}(BO(q);\mathbb{F}_{2})$ for the
classifying spaces of finite dimensional vector bundles \cite{milnor}, and in
particular (Theorem \ref{T:finitebo}) we analyze $H^{\ast}(BO(2^{n+1}%
-1);\mathbb{F}_{2})$.

Finally, we shall produce an $\mathcal{A}$-algebra epimorphism from
$\mathcal{S}=H^{\ast}(BO;\mathbb{F}_{2})$ to each of the mod two Dickson
algebras (Theorem \ref{T:dickson}), which we characterized in \cite{ppw} as
unstable $\mathcal{A}$-algebras. In fact we shall show that the $\left(
n+1\right)  $-st Dickson algebra has the role of capturing precisely the
quotient of $\mathcal{S}=H^{\ast}(BO;\mathbb{F}_{2})$ common to the cohomology
of the $n$-th distinct connected cover $BO\left\langle \phi(n)\right\rangle $
and to $BO(2^{n+1}-1)$. We speculate about how this phenomenon may relate to
spaces beyond the range in which Dickson algebras are directly realizable topologically.

Our minimal $\mathcal{A}$-algebra presentations for all the above objects will
devolve naturally from our main presentation of $\mathcal{S}$, and in that
sense these $\mathcal{A}$-algebras are all \textquotedblleft
parallel\textquotedblright\ to the main presentation.

In Appendix I, which is independent of the rest of the paper, we
present a related result, in which the unstable $\mathcal{A}$-modules
$\mathcal{F}\left( 2^{p}-1\right)
/\mathcal{A}\overline{\mathcal{A}}_{p-2}$ appear as the filtered
quotients of a simple minimal $\mathcal{A}$-presentation for
$H^{\ast}(RP^{\infty};\mathbb{F}_{2})$. We thank Don Davis, Kathryn
Lesh, and Haynes Miller for useful conversations regarding these
modules.  We also thank John Greenlees for a stimulating conversation
leading to Remark 2.4.

The first author dedicates this paper to his parents, Daphne M. and
Eric T. Pengelley, in memoriam.

\section{Motivation, first steps, and a plan}

The unstable $\mathcal{A}$-algebra of symmetric invariants $\mathcal{S}%
=H^{\ast}(BO;\mathbb{F}_{2})$ is a polynomial algebra $\mathbb{F}_{2}%
[w_{m}\co m\geq0,w_{0}=1]$, with each elementary symmetric function
(Stiefel-Whitney class) $w_{m}$ having degree $m$ \cite{milnor}. The action of
the Steenrod algebra is completely determined from the Wu formulas
\cite{kochman,stong,wu}%
\[
Sq^{j}w_{m}=\sum_{l=0}^{j}\binom{m-j+l-1}{l}w_{j-l}w_{m+l}%
\]
and the Cartan formula on products \cite{steenrodepstein}.

To ease into our categorical point of view, and to illustrate our approach and
methods, let us begin by seeing that abstract Stiefel-Whitney classes, taken
all together as free unstable $\mathcal{A}$-algebra generators, along with
imposed \textquotedblleft Wu formulas\textquotedblright, actually
\textquotedblleft present\textquotedblright\ $\mathcal{S}$. This is something
one might easily take for granted, but should actually prove, since in
principle there might be \textquotedblleft other\textquotedblright\ relations
lurking in $\mathcal{S}$ beyond those inherent in the Wu formulas. To avoid
confusion from notational abuse, we build from abstract classes $t_{m}$ which
will correspond to the actual Stiefel-Whitney classes under an isomorphism.

\begin{sloppypar}
\begin{proposition}
[Wu formulas present $\mathcal{S}$]\label{C:BOdetermined}The unstable
$\mathcal{A}$-algebra $\mathcal{S}=H^{\ast}(BO;\mathbb{F}_{2})$ is isomorphic
to the quotient of the abstract free unstable $\mathcal{A}$-algebra on classes
$t_{m}$ in each degree $m\geq1$, modulo the left $\mathcal{A}$-ideal generated
by abstract ``Wu formulas''\ formed by writing $t$'s in place of $w$'s in the
Wu formulas above.
\end{proposition}
\end{sloppypar}

\begin{proof}
Iterating the abstract Wu formulas via the Cartan formula shows that the
abstract classes $\left\{  t_{m}\co m\geq1\right\}  $ actually generate the
abstract $\mathcal{A}$-algebra quotient considered merely as an algebra, i.e.,
its (algebra) indecomposable quotient has rank at most one in each degree. On
the other hand, by its construction the abstract $\mathcal{A}$-algebra
quotient must map onto $\mathcal{S}$ by sending each $t_{m}$ to $w_{m}$, since
the respective Wu formulas correspond. Thus the two must be isomorphic, since
$\mathcal{S}$ is free as a commutative algebra.
\end{proof}

Notice, however, that this presentation of $\mathcal{S}$ is far from minimal
in the category of unstable $\mathcal{A}$-algebras, since it used vastly more
generators than needed. What we seek instead is to achieve three features for
a minimal presentation:

{\bf Step 1}\qua Find a minimal $\mathcal{A}$-submodule of $\mathcal{S}$ that
will generate $\mathcal{S}$ as an $\mathcal{A}$-algebra.

{\bf Step 2}\qua Find a minimal presentation of this $\mathcal{A}$-submodule,
i.e., with minimal generators and minimal relations.

{\bf Step 3}\qua Form the free unstable $\mathcal{A}$-algebra $\mathcal{U}$ on
this module, and find minimal relations on $\mathcal{U}$ so that its
$\mathcal{A}$-algebra quotient produces $\mathcal{S}$.

To begin, let us find a minimal set of $\mathcal{A}$-algebra generators for
$\mathcal{S}$. Consider the (algebra) indecomposable quotient $Q\mathcal{S}$,
i.e., the vector space with basis $\{w_{m}\co m\geq1\}$ and induced $\mathcal{A}%
$-action
\[
Sq^{j}w_{m}=\binom{m-1}{j}w_{m+j}.
\]
Since $\binom{m-1}{j}$ is always zero mod two when $m+j$ is a two-power, and
never zero when $m$ is a two-power and $j$ is less than $m$, we see that the
$\mathcal{A}$-module indecomposables of $Q\mathcal{S}$ have basis exactly
$\left\{  w_{2^{k}}\co k\geq0\right\}  $.

Since our philosophy is to begin the presentation at the $\mathcal{A}$-module
level, with minimal $\mathcal{A}$-algebra generators and minimal module
relations, we thus start with

\begin{definition}
Let $\mathcal{M}$ be the free unstable $\mathcal{A}$-module on abstract
classes $\{t_{2^{k}}\co k\geq0\}$, where subscripts indicate the topological
degree of each class.
\end{definition}

We wish to map $\mathcal{M}$ to $\mathcal{S}$ via $t_{2^{k}}\rightarrow
w_{2^{k}}$, and need first to ask whether $\mathcal{M}$ injects. In other
words, is the $\mathcal{A}$-submodule of $\mathcal{S}=H^{\ast}(BO;\mathbb{F}%
_{2})$ generated by $\{w_{2^{k}}\co k\geq0\}$ free? Or are there, to the
contrary, $\mathcal{A}$-relations amongst the two-power Stiefel-Whitney
classes, which will compel us to introduce module relations on $\mathcal{M}$
in order to complete steps 1 and 2 above? The Wu formulas appear to suggest
that no such relations exist. In fact we can prove something even stronger.

\begin{sloppypar}
\begin{theorem}
[Stiefel-Whitney classes inject freely]\label{T:injects}The $\mathcal{A}%
$-submodule of $\mathcal{S}=H^{\ast}(BO;\mathbb{F}_{2})$ generated by
$\{w_{m}\co m\geq0\}$ is free unstable on these classes.
\end{theorem}
\end{sloppypar}

The proof is in Section \ref{S:proofs}.

\begin{remark}
The proof also shows that in
\[
H^{\ast}(BO(q);\mathbb{F}_{2})\cong H^{\ast}(BO;\mathbb{F}_{2})/\left(
w_{m}\co m>q\right)  ,\text{ }%
\]
the $\mathcal{A}$-submodule generated by $\{w_{m}\co 0\leq m\leq q\}$ is free
unstable on these classes.
\end{remark}

\begin{remark}
The fact that the free unstable $\mathcal{A}$-module $\mathcal{F}_{m}$ on a
single class in degree $m$ injects into $H^{\ast}(BO(m);\mathbb{F}_{2})$ on
the class $w_{m}$ is clear from the already known result \cite[page
55]{lannes-zarati} that $\mathcal{F}_{m}$ is isomorphic to the invariants
$\left(  \mathcal{F}_{1}^{\otimes m}\right)  ^{\Sigma_{m}},$ which clearly
inject naturally into $\left(  H^{\ast}(RP^{\infty};\mathbb{F}_{2})^{\otimes
m}\right)  ^{\Sigma_{m}}\cong H^{\ast}(BO(m);\mathbb{F}_{2})$ on $w_{m}.$
Theorem \ref{T:injects} generalizes this by handling all $\mathcal{F}_{m}$
simultaneously, showing that they do not interfere when simultaneously perched
on the Stiefel-Whitney classes in the symmetric algebra $\mathcal{S}=H^{\ast
}(BO;\mathbb{F}_{2}).$
\end{remark}

\begin{corollary}
The $\mathcal{A}$-submodule of $\mathcal{S}=H^{\ast}(BO;\mathbb{F}_{2})$
generated by $\{w_{2^{k}}\co k\geq0\}$ is free unstable, so $\mathcal{M}$ injects
naturally into $\mathcal{S}$.
\end{corollary}

This completes steps 1 and 2 of our goal, and we can begin step 3.

\begin{definition}
Let $\mathcal{U}$ be the free unstable $\mathcal{A}$-algebra on $\mathcal{M}$,
in other words, $\mathcal{U}$ is the free unstable $\mathcal{A}$-algebra on
abstract classes $\{t_{2^{k}}\co k\geq0\}$.
\end{definition}

Clearly $\mathcal{U}$ maps via $t_{2^{k}}\rightarrow w_{2^{k}}$ onto the
desired $\mathcal{A}$-algebra $\mathcal{S}$, but the map has an enormous
kernel, since $Q\mathcal{S}$ is the vector space $\mathbb{F}_{2}\{w_{m}%
\co m\geq1\}$, while $Q\mathcal{U}$ is much larger. Our goal in step 3 is to
describe a minimal set of $\mathcal{A}$-algebra relations producing
$\mathcal{S}$ from $\mathcal{U}$, i.e., minimal generators for the kernel as
an $\mathcal{A}$-ideal.

Let us explore a prototype example in degree five, which is the first place a
difference occurs. There $Q\mathcal{S}$ has only $w_{5}$, whereas $Sq^{1}%
t_{4}$ and $Sq^{2}Sq^{1}t_{2}$ are distinct indecomposables in $Q\mathcal{U}$
(recall that $Q\mathcal{U}\cong\Sigma\Omega\mathcal{M}$, and that a basis for
$\mathcal{M}$ consists of the unstable admissible monomials on the
$\mathcal{A}$-generators $t_{2^{k}}$ \cite{steenrodepstein}). A few
calculations with the Wu formulas show that in $\mathcal{S}$ we have
\begin{align*}
Sq^{1}w_{4}  &  =w_{5}+w_{1}w_{4}\text{ and }\\
Sq^{2}Sq^{1}w_{2}  &  =w_{5}+w_{1}w_{4}+w_{2}w_{3}+w_{1}w_{2}^{2}+w_{1}%
^{2}w_{3}+w_{1}^{3}w_{2}.
\end{align*}
Thus to imitate $\mathcal{S}$ abstractly via $\mathcal{U}$, we must impose an
algebra relation on $\mathcal{U}$ decreeing that
\[
Sq^{1}t_{4}=Sq^{2}Sq^{1}t_{2}+\,\text{some decomposable,}%
\]
per the calculations above. One challenge in doing even this, though, is that
it is not clear how to describe that needed decomposable difference in
$\mathcal{U}$, since there we have no name as yet for the element
corresponding to $w_{3}$. To remedy this, and to describe general formulas for
relationships like the one we have just discovered, we wish to use the Wu
formulas to focus our understanding as much as possible on both two-power
Steenrod squares and two-power Stiefel-Whitney classes. Thus one of our
formulas in the next section will express each Stiefel-Whitney class purely in
this way (Lemma \ref{L:magic}).

While the plethora of algebra relations, such as the one above, needed to
obtain $\mathcal{S}$ from $\mathcal{U}$ may appear intractable to specify,
recall that our chosen task is actually somewhat different. Since we are
working in the category of $\mathcal{A}$-algebras, we seek relations in
$\mathcal{U}$ whose $\mathcal{A}$-algebra consequences, not just their algebra
consequences, will produce $\mathcal{S}$. We shall show that this requires
only a much smaller and more tractable set of relations, for which our
illustration in degree five serves as perfect prototype. Specifically, the
relationship between $Sq^{2^{i}}w_{2^{k}}$ and $Sq^{2^{k-1}}Sq^{2^{i}%
}w_{2^{k-1}}$ for every $i\leq k-2$ will be the key place to focus attention.
We shall impose one abstract relation on $\mathcal{U}$ for each such pair
$(k,i)$, and prove that these are precisely the minimal relations producing
$\mathcal{S}=H^{\ast}(BO;\mathbb{F}_{2})$ in the category of $\mathcal{A}$-algebras.

Our general plan is as follows. Form our abstract presentation candidate as
just outlined; call it $\mathcal{G}$. The construction of $\mathcal{G}$ will
immediately provide a natural $\mathcal{A}$-algebra epimorphism to
$\mathcal{S}$. The hard part now is showing that our $(k,i)$-indexed family of
$\mathcal{A}$-algebra relations leaves no remaining kernel, i.e., that we have
put in enough relations to generate the kernel as an $\mathcal{A}$-ideal. To
achieve this we show that the epimorphism $\mathcal{G}\rightarrow\mathcal{S}$
induces a monomorphism $Q\mathcal{G}\rightarrow Q\mathcal{S}$, on the
indecomposable quotients, by computing a basis for $Q\mathcal{G}$. For this we
appeal to our earlier understanding \cite{ppw}, via the Kudo-Araki-May algebra
$\mathcal{K}$ \cite{pw} (see Appendix II), of bases for the unstable cyclic
$\mathcal{A}$-modules arising in the analogous structure theorem for the
Dickson algebras. With $Q\mathcal{G}\rightarrow Q\mathcal{S}$ an isomorphism,
$\mathcal{G}\rightarrow\mathcal{S}$ must be an isomorphism also, since
$\mathcal{S}$ is a free commutative algebra. The minimality of the
$(k,i)$-family of relations is then not hard to see by appropriate filtering.

\section{Main theorem}

We first identify the key $\mathcal{A}$-algebra relations in $\mathcal{S}%
=H^{\ast}(BO;\mathbb{F}_{2})$.

Analysis of the binomial coefficients in the Wu formulas shows that if
$r\geq1,$ then
\begin{equation}
Sq^{2^{j-1}}w_{r2^{j}}=w_{2^{j-1}}w_{r2^{j}}+w_{2^{j-1}+r2^{j}}.
\label{premagic}%
\end{equation}
This formula will serve two purposes. It will guide us below in how to specify
any Stiefel-Whitney class from just the two-power ones, which is needed for
creating our abstract presentation. But before this it will lead us to the key
relations needed from $\mathcal{S}$.

To find these, recall from the previous section that we seek a relation
involving a decomposable difference between $Sq^{2^{i}}w_{2^{k}}$ and
$Sq^{2^{k-1}}Sq^{2^{i}}w_{2^{k-1}}$ for every $i\leq k-2$. We begin with a
special case of equation (\ref{premagic}): For $i\leq k-2,$ we have
\[
Sq^{2^{i}}w_{2^{k-1}}=w_{2^{i}}w_{2^{k-1}}+w_{2^{k-1}+2^{i}}.
\]
Applying $Sq^{2^{k-1}},$ we get
\[
Sq^{2^{k-1}}Sq^{2^{i}}w_{2^{k-1}}=Sq^{2^{k-1}}\left(  w_{2^{i}}w_{2^{k-1}%
}\right)  +Sq^{2^{k-1}}\left(  w_{2^{k-1}+2^{i}}\right)  \text{.}%
\]
Using a Wu formula on the last term, analyzing the binomial coefficients, and
using (\ref{premagic}) again, the reader may check that we obtain the
following relations.

\begin{proposition}
[Key relations in $\mathcal{S}$]For $i\leq k-2$,%
\begin{multline}
Sq^{2^{k-1}}Sq^{2^{i}}w_{2^{k-1}}=Sq^{2^{i}}w_{2^{k}}+\label{equivmoddec}\\
Sq^{2^{k-1}}\left(  w_{2^{i}}w_{2^{k-1}}\right)  +\sum_{l=0}^{2^{k-i-1}%
-2}w_{2^{k-1}-2^{i}l}w_{2^{k-1}+2^{i}+2^{i}l},
\end{multline}
\end{proposition}

These show explicitly how the elements $Sq^{2^{i}}w_{2^{k}}$ and $Sq^{2^{k-1}%
}Sq^{2^{i}}w_{2^{k-1}}$ differ by a decomposable, and will guide us to the
corresponding abstract relations needed in $\mathcal{G}$. However, the
relations we have found here involve non-two-power Stiefel-Whitney classes,
which still have as yet no analogs in $\mathcal{U}$. We remedy this problem
now by extending equation (\ref{premagic}).

Mixing notations, we write (\ref{premagic}) as
\[
w_{2^{j-1}+r2^{j}}=(Sq^{2^{j-1}}+w_{2^{j-1}})w_{r2^{j}}%
\]
(i.e., $\left(  Sq^{m}+w_{m}\right)  x$ means \ $Sq^{m}x+w_{m}\cdot x$). The
following lemma is then immediate.

\begin{lemma}
[Expressing Stiefel-Whitney classes]\label{L:magic}Every Stiefel-Whitney class
can be expressed in terms of two-power classes and two-power squares as
follows: If we write any $m=2^{n_{1}}+\cdots+2^{n_{s}},$ where $n_{1}%
>\cdots>n_{s}$, we have
\begin{equation}
w_{m}=\left(  Sq^{2^{n_{s}}}+w_{2^{n_{s}}}\right)  \cdots\left(  Sq^{2^{n_{2}%
}}+w_{2^{n_{2}}}\right)  w_{2^{n_{1}}}. \label{frankmagic}%
\end{equation}
\end{lemma}

We are now ready to define formally the abstract presentation $\mathcal{G}$.

\begin{definition}
\label{D:magicclasses}In $\mathcal{U}$, extend the set of generators $\left\{
t_{2^{k}},k\geq0\right\}  $, to define elements $t_{m}$ for all $m\geq1$,\ by
first writing $m=2^{n_{1}}+\cdots+2^{n_{s}},$ where \ $n_{1}>\cdots>n_{s}$.
\ Then by analogy with equation (\ref{frankmagic}) set
\[
t_{m}=(Sq^{2^{n_{s}}}+t_{2^{n_{s}}})\cdots\left(  Sq^{2^{n_{2}}}+t_{2^{n_{2}}%
}\right)  t_{2^{n_{1}}}.
\]
\end{definition}

\begin{definition}[\textup{(Abstract key relations)}] 
Imitating equation (\ref{equivmoddec}), let
$\mathcal{G}$\ be the the $\mathcal{A}$-algebra quotient of $\mathcal{U}$ by
the left\ $\mathcal{A}$-ideal generated by the elements
\begin{multline}
\theta\left(  k,i\right)  =Sq^{2^{i}}t_{2^{k}}+Sq^{2^{k-1}}Sq^{2^{i}%
}t_{2^{k-1}}+\label{equivmoddecabstract}\\
Sq^{2^{k-1}}\left(  t_{2^{k-1}}t_{2^{i}}\right)  +\sum_{l=0}^{2^{k-i-1}%
-2}t_{2^{k-1}-2^{i}l}t_{2^{k-1}+2^{i}+2^{i}l}%
\end{multline}
for $i\leq k-2.$
\end{definition}

\begin{theorem}
[Structure of $\mathcal{S}$]\label{T:mainthm}The symmetric algebra
$\mathcal{S}=H^{\ast}(BO;\mathbb{F}_{2})$\ is isomorphic to $\mathcal{G}$\ as
an algebra over the Steenrod algebra. Moreover, the relations
(\ref{equivmoddecabstract}) generating the $\mathcal{A}$-ideal are minimal,
i.e., nonredundant.
\end{theorem}

The proof is in Section \ref{S:proofs}.

\section{Applications and speculation}

We apply the main structure theorem to the cohomology images from the
connected covers of $BO$, and to the cohomology of the spaces $BO(q)$ for
classifying finite dimensional vector bundles. Finally we shall see how these
descriptions naturally converge into the Dickson invariant algebras.

First we consider cohomology images from the connected covers.

\begin{definition}
Following \cite{kochman}, let $B^{\ast}(n)$ be the cohomology image of the map
induced by the projection
\[
BO\left\langle \phi(n)\right\rangle \rightarrow BO,
\]
where $BO\!\left\langle \phi(n)\right\rangle $ is the $n$-th distinct connected
cover of $BO$. That is, $BO\!\left\langle \phi(n)\right\rangle $ is $\left(
\phi(n)-1\right)  $-connected, where $n=4s+t$, $0\leq t\leq3$, and
$\phi(n)=8s+2^{t}$.
\end{definition}

In particular, for $n=0,1,2,3$ the projections are surjective in cohomology,
so the unstable $\mathcal{A}$-algebras $B^{\ast}(n)$ are isomorphic to the
cohomologies of $BO$, $BSO$, $BSpin$, and $BO\left\langle 8\right\rangle $
\cite{kochman}. In general, $B^{\ast}(n)$ is $(2^{n}-1)$-connected, and is the
quotient of $B^{\ast}(0)=H^{\ast}BO=\mathcal{S}$ by the $\mathcal{A}$-ideal
generated by $\{w_{2^{k}}\co k<n\}$ \cite{kochman}.

\begin{theorem}
[Structure of connected cover images]\label{T:conncover}An abstract
presentation of $B^{\ast}(n)$ is obtained from that of $B^{\ast}(0)=H^{\ast
}BO=\mathcal{S}$ (Theorem \ref{T:mainthm}) as the quotient by the
$\mathcal{A}$-ideal generated by $\{t_{2^{k}}\co k<n\}$. This produces a minimal
presentation as follows.

Let $K_{n}$ denote the direct sum of the $\mathcal{A}$-module $\mathcal{M}%
(n,0)$ on $t_{2^{n}}$ with the free unstable $\mathcal{A}$-module on the
$t_{2^{k}},$ $k\geq n+1.$ \ Here $\mathcal{M}(n,0)$ is as defined in
\cite{ppw}, namely the free unstable $\mathcal{A}$-module on one generator
$t_{2^{n}}$ modulo the left $\mathcal{A}$-submodule generated by $Sq^{2^{i}%
}t_{2^{n}},$ $i\leq n-2.$

Then $B^{\ast}(n)$ is isomorphic to the quotient of the free unstable
$\mathcal{A}$-algebra on $K_{n}$ by the left\ $\mathcal{A}$-ideal generated by
the elements $\theta\left(  k,i\right)  $, $k\geq n+1$, $i\leq k-2$, subject
to the requirement that all appearances in $\theta\left(  k,i\right)  $ of
$t_{m},$ $0<m<2^{n},$ are replaced by zero.
\end{theorem}

The proof is in Section \ref{S:proofs}.

For our second application, we note that the presentation for $H^{\ast}BO$ in
our main theorem will immediately produce presentations for the cohomologies
of the classifying spaces $H^{\ast}BO(q)$, since each is just the algebra
quotient (actually also $\mathcal{A}$-algebra quotient) of $H^{\ast}BO$ by the
ideal generated by $\{w_{m}\co m>q\}$ \cite{milnor}, and $w_{m}$ corresponds to
$t_{m}$, which we defined in the presentation of $H^{\ast}BO$. The resulting
presentation becomes both tractable and useful for $H^{\ast}BO(2^{n+1}-1)$.

\begin{theorem}
[Structure of $H^{\ast}BO(2^{n+1}-1)$]\label{T:finitebo}An abstract
presentation of $H^{\ast}BO(2^{n+1}-1)$ is obtained from that of $B^{\ast
}(0)=H^{\ast}BO=\mathcal{S}$ (Theorem \ref{T:mainthm}) as the quotient by the
$\mathcal{A}$-ideal generated by $\{t_{2^{k}}\co k\geq n+1\}$. This produces a
minimal presentation as follows.

$H^{\ast}BO(2^{n+1}-1)$ is presented by the free unstable $\mathcal{A}%
$-algebra on abstract classes $\{t_{2^{k}}\co 0\leq k\leq n\}$, modulo the
left\ $\mathcal{A}$-ideal generated by the elements $\theta\left(  k,i\right)
$ for $k\leq n+1$, $i\leq k-2$, (using Definition \ref{D:magicclasses} of
$t_{m}$ for $m<2^{n+1}$), subject to the requirement that when $k=n+1$, the
term $Sq^{2^{i}}t_{2^{n+1}}$ is replaced by zero for each $i$ (all other terms
involve only $t$'s in degrees less than $2^{n+1}$).
\end{theorem}

The proof is in Section \ref{S:proofs}.

Finally, combining the relations on $\mathcal{S}=H^{\ast}(BO;\mathbb{F}_{2})$
from the two theorems above will produce the common $\mathcal{A}$-algebra
quotient of $B^{\ast}(n)$ and $H^{\ast}BO(2^{n+1}-1)$. Since the first of
these is $\left(  2^{n}-1\right)  $-connected, while the second is
decomposable beyond degree $2^{n+1}-1$, we will obtain an $\mathcal{A}%
$-algebra with algebra generators in the range $2^{n}$ through $2^{n+1}-1.$
Surprisingly, this much smaller quotient of $\mathcal{S}=H^{\ast}BO$ turns out
to be already familiar. We will show now that as an $\mathcal{A}$-algebra it
is isomorphic to the $n$-th Dickson algebra $W_{n+1}$ (see Figure
\ref{F:One}). In this sense one can say that the Dickson algebra captures
precisely the cohomology common to $BO\left\langle \phi(n)\right\rangle $ and
$BO(2^{n+1}-1)$ from $H^{\ast}BO$, i.e., it is the $\mathcal{A}$-algebra pushout.

\begin{figure}[tbh]
\[%
\begin{array}
[c]{ccc}%
W_{n+1} & \leftarrow & H^{\ast}BO(2^{n+1}-1)\\
&  & \\
\uparrow &  & \uparrow\\
&  & \\
B^{\ast}(n) & \leftarrow & H^{\ast}BO
\end{array}
\]\caption{}%
\label{F:One}%
\end{figure}

\begin{theorem}
[Convergence to Dickson algebras]\label{T:dickson}The quotient of the
symmetric algebra $\mathcal{S}$ by the left\ $\mathcal{A}$-ideal generated by
$\left\{  w_{2^{k}}\co k\neq n\right\}  $ is isomorphic to the $n+1$-st mod $2$
Dickson algebra, $W_{n+1}$. Specifically, using the notation of the
presentation of Theorem \ref{T:mainthm}, as an $\mathcal{A}$-algebra it is
minimally presented by the free unstable $\mathcal{A}$-algebra on the module
$\mathcal{M}(n,0)$ (defined in Theorem \ref{T:conncover}), subject to the
single $\mathcal{A}$-algebra relation
\[
Sq^{2^{n}}Sq^{2^{n-1}}t_{2^{n}}=t_{2^{n}}Sq^{2^{n-1}}t_{2^{n}}.
\]
We proved in \cite{ppw} that this precisely characterizes the Dickson algebra
$W_{n+1}$.
\end{theorem}

The proof is in Section \ref{S:proofs}.

Let us speculate on how Figure \ref{F:One} might fit in with something
topologically realizable. It is known that $W_{n+1}$ is realizable precisely
for $n\leq3$ \cite{realize}, and that $B^{\ast}(n)\subset H^{\ast
}BO\left\langle \phi(n)\right\rangle $ is an isomorphism also precisely in
this range \cite{kochman}. Thus for $n\leq3$ it is reasonable to expect that
Figure \ref{F:One} be realizable. For general $n$ it is perhaps reasonable to
hope for the existence of a space $X_{n}$ and a homotopy commutative square
(Figure \ref{F:Two}) whose cohomology is compatible with Figure \ref{F:One} in
the sense of combining to produce the commutative diagram of Figure
\ref{F:Three}. Additionally we would like $X_{n}$ to have the property that
the outer square in Figure \ref{F:Three} is also a pushout of unstable
$\mathcal{A}$-algebras. In other words, $X_{n}$ does its best to realize a
Dickson algebra, even when this is no longer possible.

\begin{figure}[tbh]
\[%
\begin{array}
[c]{ccc}%
X_{n} & \rightarrow & BO(2^{n+1}-1)\\
&  & \\
\downarrow &  & \downarrow\\
&  & \\
BO\left\langle \phi(n)\right\rangle  & \rightarrow & BO
\end{array}
\]
\caption{{}}%
\label{F:Two}%
\end{figure}

\begin{figure}[tbh]
\[%
\begin{array}
[c]{ccccc}%
H^{\ast}X_{n} & \leftarrow & W_{n+1} & \leftarrow & H^{\ast}BO(2^{n+1}-1)\\
&  &  &  & \\
\uparrow &  & \uparrow &  & \uparrow\\
&  &  &  & \\
H^{\ast}BO\left\langle \phi(n)\right\rangle  & \supset & B^{\ast}(n) &
\leftarrow & H^{\ast}BO
\end{array}
\]
\caption{{}}%
\label{F:Three}%
\end{figure}

\section{Proofs\label{S:proofs}}

\begin{proof}
[Proof of Theorem \ref{T:injects}]Let $\mathcal{F}_{m}$ be the free unstable
$\mathcal{A}$-module (equivalently $\mathcal{K}$ module) on a generator
$t_{m}$ in degree $m$. We shall show that the $\mathcal{A}$-module map
$f\co \oplus_{m\geq0}\mathcal{F}_{m}\rightarrow H^{\ast}BO$ determined by
$f(t_{m})=w_{m}$ is injective.

{}From \cite{pw}, basis elements for the domain of $f$ consist of $D_{J}t_{m}$
where $J=\left(  j_{1},\ldots,j_{s}\right)  $ and $0\leq j_{1}\leq\cdots\leq
j_{s}<m$. (Appendix II recalls the features of  the elements $D_{J}$ in the
Kudo-Araki-May algebra $\mathcal{K}$ essential to what follows.)

On the other side of $f$, basis monomials of the range $H^{\ast}BO$ can be
written as $\cdots w_{n_{2}}w_{n_{1}}$ with nondecreasing indices, i.e.,
labeled by finitely nonzero tuples $\left(  \ldots,n_{2},n_{1}\right)  $ with
$0\leq\cdots\leq n_{2}\leq n_{1}$. We order the latter reverse lexicographically.

Now for each basis element $D_{J}t_{m}$, we consider its image $f\left(
D_{J}t_{m}\right)  =D_{J}w_{m}$, and we claim that this element of $H^{\ast
}BO$ has a \textquotedblleft leading\textquotedblright\ monomial term, i.e.,
that
\[
D_{J}w_{m}=\underset{z}{\underbrace{w_{m-j_{s}}^{2^{s-1}}w_{m-j_{s-1}%
}^{2^{s-2}}\cdots w_{m-j_{2}}^{2}w_{m-j_{1}}w_{m}}}+\text{ higher order
terms.}%
\]
This will complete the proof, since distinct $D_{J}w_{m}$ clearly produce
distinct leading monomials, with remaining terms always of higher order; so
the $D_{J}w_{m}$ are all linearly independent, and thus $f$ is injective.

We will use the following notation: As a subscript, ``$>k$''\ (resp.\ ``$<k$'')
denotes any index greater (resp.\ less) than $k$, each occurrence of an
unsubscripted $w$ denotes any element of $H^{\ast}BO$, and expressions
involving any of these mean any sum of expressions of such form.

We prove our claim by induction on $s$, based on the Wu formula
\[
D_{j}w_{m}=Sq^{m-j}w_{m}=w_{m-j}w_{m}\text{ }+\text{ higher order terms of
form }ww_{>m}\text{.}%
\]
Clearly the claim holds for lengths $0$ and $1$. For the inductive step,
consider $D_{\widehat{J}}$ of length $s+1$, and note that application of any
nontrivially-acting $D_{J}$ always increases the order of a monomial in
$H^{\ast}BO$. Now calculate, using the $\mathcal{K}$-Cartan formula \cite{pw}
as needed, and recalling that the leading term $z$ was defined above:%
\begin{align*}
D_{\widehat{J}}w_{m} &  \\
&  =D_{j_{1}}D_{j_{2}}\cdots D_{j_{s+1}}w_{m}=D_{j_{1}}\left(  D_{j_{2}}\cdots
D_{j_{s+1}}w_{m}\right)  \\
&  =D_{j_{1}}\left(  \underset{x}{\underbrace{w_{m-j_{s+1}}^{2^{s-1}}\cdots
w_{m-j_{2}}}}w_{m}+\text{ higher order terms than }xw_{m}\right)  \\
&  =x^{2}w_{m-j_{1}}w_{m}+x^{2}ww_{>m}+wD_{<j_{1}}w_{m}\\
&  +D_{j_{1}}\left(  ww_{>m}+\text{ }\underset{v}{\underbrace{\text{higher
order terms than }x}}w_{m}\right)  \\
&  =z+ww_{>\left(  m-j_{1}\right)  }w_{m}+\text{ higher order terms than }z\\
&  +\left(  ww_{>m}+v^{2}w_{m-j_{1}}w_{m}+ww_{>\left(  m-j_{1}\right)  }%
w_{m}\right)  \\
&  =z+\text{ higher order terms than }z\text{, }%
\end{align*}
since the terms of $v^{2}$ have higher order than $x^{2}$.
\end{proof}

\begin{proof}
[Proof of Theorem \ref{T:mainthm}]There is a map of $\mathcal{A}$-algebras
$\mathcal{U\rightarrow S}$ obtained by taking $t_{2^{k}}$ to\ $w_{2^{k}},$ and
from Lemma \ref{L:magic} and Definition \ref{D:magicclasses} this map takes
each $t_{m}$ to $w_{m}.$ \ Since the relations (\ref{equivmoddecabstract})
that define $\mathcal{G}$ map to those also satisfied in $\mathcal{S}$
(\ref{equivmoddec}), there is an induced $\mathcal{A}$-algebra epimorphism
$\mathcal{G}\rightarrow\mathcal{S}$. We shall show that this map is monic by
showing that the induced map on the indecomposable quotients is monic,
essentially a counting argument.

To start with, note that the indecomposables are
\[
Q\mathcal{U=}\left\langle Sq^{I}t_{2^{k}}\co k\geq0,\text{ }I\text{ admissible,
of excess }<2^{k}\right\rangle .
\]
Then $Q\mathcal{G}$\ is\ $Q\mathcal{U}$\ modulo the $\mathcal{A}$-relations
(degenerate versions of $\theta(k,i)=0$)%
\[
Sq^{2^{i}}t_{2^{k}}=Sq^{2^{k-1}}Sq^{2^{i}}t_{2^{k-1}},\ i\leq k-2.\
\]
There is an\ $\mathcal{A}$-module filtration \
\[
F_{p}Q\mathcal{U=}\left\langle Sq^{I}t_{2^{k}}\co 0\leq k\leq p\text{, }I\text{
admissible, of excess }<2^{k}\right\rangle ,
\]
which induces an $\mathcal{A}$-module filtration $F_{p}Q\mathcal{G}$. \ Then
\begin{multline*}
F_{p}Q\mathcal{G}/F_{p-1}Q\mathcal{G}=\\
\left\langle Sq^{I}t_{2^{p}}\co I\text{ admissible, of excess }<2^{p}%
\right\rangle /\mathcal{A}\left\{  Sq^{2^{i}}t_{2^{p}}\co i\leq p-2\right\}  .
\end{multline*}
This is the suspension of the module $\mathcal{M}(p,1)$ analyzed in
\cite[Theorem 2.11]{ppw}\footnote{$\mathcal{M}(p,1)$ is defined in \cite{ppw}
as the quotient of the free unstable $\mathcal{A}$-module on a class in degree
$2^{p}-1$ modulo the action of $Sq^{2^{i}}$ for $i\leq p-2$; in other words,
in usual notation, $\mathcal{M}(p,1)=\mathcal{F}\left(  2^{p}-1\right)
/\mathcal{A}\overline{\mathcal{A}}_{p-2}.$}, and the basis described there
suspends to
\[
\left\{  D_{I}t_{2^{p}}\co I=(2^{a_{1}},\ldots,2^{a_{l}}),\text{ where }0\leq
a_{1}\leq\cdots\leq a_{l}<p\right\}  .
\]
(As in the proof of Theorem \ref{T:injects}, we refer the reader to Appendix
II for essentials concerning the elements $D_{I}$ in the Kudo-Araki-May
algebra $\mathcal{K}$.)

We shall finish the proof of isomorphism by showing that the above basis
elements for $\oplus_{p\geq0}F_{p}Q\mathcal{G}/F_{p-1}Q\mathcal{G}$ are in
distinct degrees; in fact we claim there is exactly one in each positive
degree (The appendix discusses the modules $\mathcal{M}(p,1)$ in relation to
the literature, and points out an alternative path for substantiating our
claim.).\ Let $m$ be a positive integer. \ Then $m$ may be written uniquely in
the form
\[
m=2^{r}-\sum_{j=1}^{s}2^{b_{j}},
\]
where $s\geq0$ and $0\leq b_{1}<\cdots<b_{s}<r-1$. \ The reader may check by
induction on $s$ that the unique basis element in degree $m$ is $D_{I}%
t_{2^{p}},$ where $p=r-s$ and $I=(2^{a_{1}},\ldots,2^{a_{s}}),$ with
$a_{j}=b_{j}-j+1$. With both $Q\mathcal{G}$ and $Q\mathcal{S}$ having rank one
in each degree, $Q\mathcal{G}\rightarrow Q\mathcal{S}$ is an isomorphism. Then
since $\mathcal{S}$ is a free commutative algebra, the epimorphism
$\mathcal{G}\rightarrow\mathcal{S}$ must be an isomorphism also.

That the relations are minimal (nonredundant) is clear from the fact that in
$F_{p}Q\mathcal{U}/F_{p-1}Q\mathcal{U}$, which is the suspension of the free
unstable module on a class in degree $2^{p}-1$, the induced relations are
simply $Sq^{2^{i}}t_{2^{p}}=0$, for $i\leq p-2$, and these are all nonredundant.
\end{proof}

\begin{proof}
[Proof of Theorem \ref{T:conncover}]We have already mentioned that according
to \cite{kochman}, $B^{\ast}(n)$ is isomorphic to the quotient of
$\mathcal{S}$\ by the $\mathcal{A}$-ideal generated by $\{w_{2^{k}}\co k\leq
n-1\}.$ \ Hence the images under the projection $\mathcal{S\rightarrow}%
B^{\ast}(n)$ of all $w_{m},1\leq m\leq2^{n}-1,$ are certainly zero from Lemma
\ref{L:magic}. \ {}From \cite{kochman} we also have that $B^{\ast}(n)$ is a
polynomial algebra generated by certain remaining $w_{m}$ (see below). We
denote the images of the $w_{m}$ in $B^{\ast}(n)$ by the same symbols $w_{m}.$ \ 

Let $\mathcal{H}_{n}$ denote the quotient of the free unstable $\mathcal{A}%
$-algebra on $K_{n}$ by the left\ $\mathcal{A}$-ideal generated by the
elements $\theta\left(  k,i\right)  $, for $k\geq n+1$, subject to the
requirement that all appearances of $t_{m},$ $0<m<2^{n}$, are replaced by
zero, as in the statement of the theorem.

We begin by defining a map from $K_{n}$ to $B^{\ast}(n)$ by, as in the
preceding proof, assigning $t_{2^{k}}$ to $w_{2^{k}}$ for $k\geq n.$ \ Since
the defining relations for $K_{n}$ are clearly satisfied in $B^{\ast}(n)$
(from equation (\ref{equivmoddec})), this assignment extends to the desired
map.  And since the defining relations for the algebra $\mathcal{H}_{n}$ are
also clearly satisfied in $B^{\ast}(n)$, this extends to an $\mathcal{A}%
$-algebra map $\mathcal{H}_{n}\rightarrow B^{\ast}(n).$ \ This map is
epimorphic (since $B^{\ast}(n)$ is generated by certain $w_{m}$ with
$i\geq2^{n}$), so as in the preceding proof, we need only show the the induced
map on indecomposables is monomorphic.  

According to \cite{kochman}\footnote{Kochman describes degrees of generators
in terms of $\alpha(m)+\nu(m)$ ($\nu$ is the $2$-divisibility), but we
equivalently use $\alpha(m-1)=\alpha(m)+\nu(m)-1$.}, the polynomial generators
of $B^{\ast}(n)$ are the $w_{m}$ for which $\alpha(m-1),$ the number of ones
in the binary representation of $m-1,$ is at least $n$. \ We filter
$Q\mathcal{H}_{n}$ as in the proof of the previous theorem,
\[
F_{p}Q\mathcal{H}_{n}\mathcal{=}\left\langle Sq^{I}t_{2^{k}}\in Q\mathcal{H}%
_{n}\co k\leq p\right\rangle ,
\]
and as in the previous proof the filtered quotient $F_{p}Q\mathcal{H}%
_{n}/F_{p-1}Q\mathcal{H}_{n}$ is the suspension of the module $\mathcal{M}%
(p,1)$ for $p\geq n$, and $0$ for $p<n$. \ It is straightforward to check that
the alpha numbers of one less than the degrees of the elements
\[
\left\{  D_{I}t_{2^{p}}\co I=(2^{a_{1}},\ldots,2^{a_{l}}),\text{ where }0\leq
a_{1}\leq\cdots\leq a_{l}<p\right\}
\]
are exactly $p\geq n$, so these are all in degrees where $B^{\ast}(n)$ has
generators. \ Since we showed in the previous proof that these elements are
also in distinct degrees, this similarly completes the proof. Minimality
follows as in the previous proof.
\end{proof}

\begin{sloppypar}
\begin{proof}
[Proof of Theorem \ref{T:finitebo}]It is clear that the presentation of
$\mathcal{S}$ collapses in the manner stated. Minimality follows for most of
the relations as in the previous proofs. We comment only that to confirm that
the collapsed top relations
\[
0=\theta\left(  n+1,i\right)  \equiv Sq^{2^{n}}Sq^{2^{i}}t_{2^{n}%
}+\text{decomposables for }i\leq n-1
\]
are also all nonredundant, one can observe that there is a natural map of the
new presentation without these final relations to the presentation for
$\mathcal{S}$, and compute that on indecomposables, each $Sq^{2^{n}}Sq^{2^{i}%
}t_{2^{n}}$ maps to $w_{2^{n+1}+2^{i}}$. Now from the Wu formulas, $QH^{\ast
}BO$ is filtered over $\mathcal{A}$ by $F_{p}QH^{\ast}BO=\left\{  w_{m}%
\co \alpha(m-1)\leq p\right\}  $, and $w_{2^{n+1}+2^{i}}$ is in filtration
exactly $i+1$. Thus $\left\{  w_{2^{n+1}+2^{i}}\co i\leq n-1\right\}  $ 
must be a
minimal generating set for the $\mathcal{A}$-submodule it generates in
$QH^{\ast}BO$. The same then must be true of $\left\{  \theta\left(
n+1,i\right)  \co i\leq n-1\right\}  $ in the indecomposables of the new
presentation without these final relations; so they too are minimal.
\end{proof}
\end{sloppypar}

\begin{proof}
[Proof of Theorem \ref{T:dickson}]In \cite{ppw} we proved that the $\left(
n+1\right)  $-st Dickson algebra $W_{n+1}$ is isomorphic to the quotient of
the free unstable $\mathcal{A}$-algebra on the module $\mathcal{M}(n,0)$ on
generator $x_{2^{n}}$ by the single $\mathcal{A}$-algebra relation
\[
Sq^{2^{n}}Sq^{2^{n-1}}x_{2^{n}}=x_{2^{n}}Sq^{2^{n-1}}x_{2^{n}},
\]
and that $\mathcal{M}(n,0)$ injects into $W_{n+1}$ (\cite{ppw}, proof of
Theorem 2.11). In other words, this is a minimal presentation in our sense.

Now let us turn to the quotient of the symmetric algebra that combines the
relations from the previous two theorems, i.e., the quotient by the
left\ $\mathcal{A}$-ideal generated by $\left\{  t_{2^{k}},k\neq n\right\}  $.
Let us denote this quotient by $\mathcal{J}_{n}.$ \ In $\mathcal{J}_{n},$ the
relations $\theta\left(  k,i\right)  $ are all trivial except when $k$ is
$n+1$ or $n.$ \ When $k=n,$ they reduce to $Sq^{2^{i}}t_{2^{n}}=0$, $i\leq
n-2,$ the defining relations for $\mathcal{M}(n,0).$ \ When $k=n+1$, we have
the relations
\begin{multline*}
0=\theta\left(  n+1,i\right)  \equiv Sq^{2^{n}}Sq^{2^{i}}t_{2^{n}}+Sq^{2^{i}%
}t_{2^{n+1}}+\\
Sq^{2^{n}}\left(  t_{2^{n}}t_{2^{i}}\right)  +\sum_{l=0}^{2^{n-i}-2}%
t_{2^{n}-2^{i}l}t_{2^{n}+2^{i}+2^{i}l}%
\end{multline*}
for $i\leq n-1.$ \ These reduce to
\[
Sq^{2^{n}}Sq^{2^{i}}t_{2^{n}}=t_{2^{n}}t_{2^{n}+2^{i}}.
\]
Now since
\[
t_{2^{n}}t_{2^{n}+2^{i}}=t_{2^{n}}\left(  Sq^{2^{i}}t_{2^{n}}+t_{2^{i}%
}t_{2^{n}}\right)  =t_{2^{n}}Sq^{2^{i}}t_{2^{n}},
\]
the relations can be rewritten as
\[
Sq^{2^{n}}Sq^{2^{i}}t_{2^{n}}=t_{2^{n}}Sq^{2^{i}}t_{2^{n}}.
\]
Since $Sq^{2^{i}}t_{2^{n}}=0$ for $i<n-1$, these are trivial for $i<n-1$, and
yield
\[
Sq^{2^{n}}Sq^{2^{n-1}}t_{2^{n}}=t_{2^{n}}Sq^{2^{n-1}}t_{2^{n}}%
\]
for $i=n-1$. This precisely matches the single relation (stated above)
characterizing the Dickson algebra, so we obtain an isomorphism of
$\mathcal{A}$-algebras from $\mathcal{J}_{n}$ to $W_{n+1}$ by taking
$t_{2^{n}}\in\mathcal{J}_{n}$ to the generator $x_{2^{n}}\in W_{n+1}.$
\end{proof}

\section{Appendix I: The unstable modules $\mathcal{F}\!\left(  2^{p}{-}1\right)
/\mathcal{A}\overline{\mathcal{A}}_{p-2}$ and a minimal $\mathcal{A}%
$-presentation for $H^{\ast}\!\left(  RP^{\infty}\right)  $}

For each $p\geq0$ , the module $\mathcal{M}(p,1)$ is defined in \cite{ppw} as
the quotient of the free unstable $\mathcal{A}$-module on a class $x_{2^{p}%
-1}$ in degree $2^{p}-1$ modulo the action of $Sq^{2^{i}}$ for $i\leq p-2$; in
other words, in usual notation,
\[
\mathcal{M}(p,1)=\mathcal{F}\left(  2^{p}-1\right)  /\mathcal{A}%
\overline{\mathcal{A}}_{p-2}.
\]
These modules are tractable, important, and interesting, and we shall show
they are the filtered quotients of a simple minimal $\mathcal{A}$-presentation
for $H^{\ast}RP^{\infty}$.

In the proof of our primary Theorem \ref{T:mainthm} above, we appealed to our
development in \cite[Theorem 2.11]{ppw} of bases for these modules. The proof
used the bases to \textquotedblleft count\textquotedblright\ that the direct
sum of the modules (we were actually dealing with their suspensions in that
theorem) has rank exactly one in each nonnegative degree. In fact we know the
rank separately for each module:

\begin{theorem}
[Rank of $\mathcal{M}(p,1)$]\label{T:modulerank}The module $\mathcal{M}(p,1)$
has precisely a single nonzero element in each degree with alpha number $p$,
i.e., with $p$ ones in its binary expansion, and nothing else.
\end{theorem}

\begin{proof}
The basis for $\mathcal{M}(p,1)$ provided in \cite[Theorem 2.11]{ppw} is
\begin{align*}
\{D_{I}x_{2^{p}-1}\co \text{ } &  \text{the multi-index }I\text{ consists of
nonnegative, }\\
&  \text{nondecreasing entries of form }2^{k}-1\text{, }k<p\}.
\end{align*}
The reader may check that the degrees of these elements are precisely those
with alpha number $p$ (see Appendix II for a recollection of essentials
regarding the elements $D_{I}$ in the Kudo-Araki-May algebra $\mathcal{K}$).
\end{proof}

This suggests a connection to the cohomology of $RP^{\infty}$. Recall that
\begin{equation}
H^{\ast}RP^{\infty}\cong\mathbb{F}_{2}[y]\text{ with }Sq^{j}y^{l}=\binom{l}%
{j}y^{l+j}, \label{projectivespace}%
\end{equation}
from which one sees that $H^{\ast}RP^{\infty}$ is $\mathcal{A}$-filtered by
the number of ones in the binary expansion of degrees. Indeed it is now not
hard to prove

\begin{theorem}
[$\mathcal{M}(p,1)$ and $H^{\ast}RP^{\infty}$]\label{T:moduleiso}The
$\mathcal{A}$-module $\mathcal{M}(p,1)$ is isomorphic to the $p$-th filtered
quotient of $H^{\ast}RP^{\infty}$.
\end{theorem}

\begin{proof}
The module $\mathcal{M}(p,1)$ clearly maps nontrivially to the $p$-th filtered
quotient of $H^{\ast}RP^{\infty}$, since the quotient begins with $y^{2^{p}%
-1}$, and $Sq^{2^{i}}y^{2^{p}-1}$ lies in lower filtration for $i\leq p-2$.
The map is onto because one sees from (\ref{projectivespace}) that the $p$-th
filtered quotient of\ $H^{\ast}RP^{\infty}$ is generated over $\mathcal{A}$
from degree $2^{p}-1$. Now the previous theorem shows that the ranks agree, so
the two are isomorphic.
\end{proof}

\begin{remark}
This result also follows from \cite{davis}, where it essentially appears in a
stabilized form. Indeed, in \cite{davis} the $\mathcal{A}$-modules
\[
\Sigma^{2^{p}-1}\mathcal{A}/\mathcal{A}\left\{  Sq^{2^{j}}\co j\neq p-1\right\}
\]
are studied with stable purposes in mind. Each of these modules obviously maps
onto the corresponding $\mathcal{M}(p,1)$, and thus the two would clearly be
isomorphic if it were known that the domain module is unstable, which does not
seem obvious. In fact, though, it is proven in \cite{davis} that these modules
are isomorphic to the same filtered quotients of $H^{\ast}RP^{\infty}$. Thus
they are indeed unstable and isomorphic to the modules $\mathcal{M}(p,1)$. The
theorem follows.
\end{remark}

\begin{remark}
The modules $\mathcal{M}(p,1)$ are also used in \cite{lesh}, where Remark 2.6
claims that in an unpublished manuscript \cite{massey}, William Massey
calculated that $\mathcal{M}(p,1)$ is $\mathcal{A}$-isomorphic to the $p$-th
filtered quotient of $H^{\ast}RP^{\infty}$, i.e., the theorem above. However,
this does not actually seem to appear explicitly in \cite{massey}. Finally, we
note that the filtered quotients of $H^{\ast}RP^{\infty}$ arise again in
\cite[after Prop.\ 3.1]{bousfield-davis} in a fashion closely related both to
\cite{lesh} and \cite{massey}.
\end{remark}

We are now equipped to show

\begin{theorem}
[Minimal $\mathcal{A}$-presentation of $H^{\ast}\left(  RP^{\infty}\right)  $%
]There is a minimal unstable $\mathcal{A}$-module presentation of $H^{\ast
}\left(  RP^{\infty};\mathbb{F}_{2}\right)  $, as the quotient of the free
unstable module on abstract classes $s_{2^{k}-1}$ in degrees $2^{k}-1$ by the
relations
\[
Sq^{2^{i}}s_{2^{k}-1}=Sq^{2^{k-1}}Sq^{2^{i}}s_{2^{k-1}-1},\ i\leq k-2.\
\]
\end{theorem}

\begin{proof}
There is an $\mathcal{A}$-module map from the abstract quotient to $H^{\ast
}RP^{\infty}$, carrying each $\mathcal{A}$-generator nontrivially, since the
given relations are easily calculated also to hold amongst the nonzero classes
in $H^{\ast}RP^{\infty}$. Moreover this is epic, since $H^{\ast}RP^{\infty}$
is generated over $\mathcal{A}$ from degrees one less than a two-power. To see
that the two are isomorphic, we need merely show that these relations are
enough, i.e., that the abstract quotient has only rank one in each degree.
This we do by considering the $\mathcal{A}$-filtration of the abstract
quotient in which the $p$-th filtration is the $\mathcal{A}$-submodule
generated by $\left\{  s_{1},\ldots,s_{2^{p}-1}\right\}  $. The $p$-th
filtered quotient is clearly $\mathcal{M}(p,1)$. That the union of these has
rank one in each nonnegative degree follows from either of the two previous theorems.

Minimality of the presentation is clear. The nonzero classes in $H^{\ast
}RP^{\infty}$ in degrees one less than a power of two cannot be reached from
below, so the generating set is minimal, and unique. The nonredundancy of all
the relations is clear from the filtered quotients and the fact that two-power
squares are minimal generators of $\mathcal{A}$.

An alternative proof would be to obtain this presentation simply by collapsing
the relations (\ref{equivmoddecabstract}) in the $\mathcal{A}$-algebra
presentation of $H^{\ast}BO$ in Theorem \ref{T:mainthm} to the indecomposable
quotient, since $\Sigma H^{\ast}RP^{\infty}\cong QH^{\ast}BO$ as $\mathcal{A}%
$-modules (Wu formulas).
\end{proof}

\section{Appendix II: The Kudo-Araki-May algebra $\mathcal{K}$}

We recall here just the bare essentials about $\mathcal{K}$ needed to
understand the proofs in this paper. We refer the reader to \cite{pw} for much
more extensive information about $\mathcal{K}$.

The mod two Kudo-Araki-May algebra $\mathcal{K}$ is the $\mathbb{F}_{2}%
$-bialgebra (with identity) generated by elements $\{D_{i}\co i\geq0\}$ subject
to homogeneous (Adem) relations \cite[Def. 2.1]{pw}, with coproduct $\phi$
determined by the formula
\[
\phi(D_{i})=\sum_{t=0}^{i}D_{t}\otimes D_{i-t}.
\]
It is bigraded by length and topological degrees ($\left\vert D_{i}\right\vert
=i$), which behave skew-additively under multiplication \cite[Def. 2.1]{pw}.

The $\mathbb{F}_{2}$-cohomology of any space is an unstable algebra over the
Steenrod algebra, and there is a correspondence between unstable $\mathcal{A}%
$-algebras and unstable $\mathcal{K}$-algebras, completely determined by
iterating the conversion formulae: On any element $x_{l}$ of degree $l$, and
for all $j\geq0$, one has
\[
D_{j}x_{l}=Sq^{l-j}x_{l}\text{, equivalently, }Sq^{j}x_{l}=D_{l-j}x_{l}.
\]
Since the degree of the element is involved in the conversion, and this
changes as operations are composed, the algebra structures of $\mathcal{A}$
and $\mathcal{K}$ are very different, and the skew additivity of the bigrading
in $\mathcal{K}$ reflects this. 

The requirements for an unstable $\mathcal{K}$-algebra, corresponding to the
nature and requirements of an unstable $\mathcal{A}$-algebra, are: On any
element $x_{l}$ of degree $l$,%
\[
D_{l}x_{l}=x_{l}\text{, }D_{j}x_{l}=0\text{ for }j>l\text{, and }D_{0}%
x_{l}=x_{l}^{2}.
\]
Finally, and used in our proofs, the $\mathcal{K}$-algebra structure obeys the
(Cartan) formula according to the coproduct $\phi$ in $\mathcal{K}$:%
\[
D_{i}(xy)=\sum_{t=0}^{i}D_{t}(x)D_{i-t}(y).
\]

\Addresses\recd 


\begin{thebibliography}{99}                                                                                               %
\bibitem {bousfield-davis}A.K. Bousfield, D.M. Davis, On the unstable Adams
spectral sequence for $SO$ and $U$, and splittings of unstable $Ext$ groups,
\emph{Boletin de la Sociedad Matematica Mexicana }(2) \textbf{37}{}(1992), 41--53.

\bibitem {davis}D. Davis, Some quotients of the Steenrod algebra, \emph{Proc.
Amer. Math. Soc.} \textbf{83} (1981), 616--618.

\bibitem {kochman}S. Kochman, An algebraic filtration of $H_{\ast}BO$, in
proc., \emph{Northwestern Homotopy Theory Conference (Evanston, Ill., 1982)},
Contemporary Math. \textbf{19} (1983), Amer. Math. Soc., 115--143.

\bibitem {lannes-zarati}J. Lannes, S. Zarati, Foncteurs d\'{e}riv\'{e}s de la
d\'{e}stabilisation, \emph{Math. Zeit.} \textbf{194} (1987), 25--59.

\bibitem {lesh}K. Lesh, A conjecture on the unstable Adams spectral sequences
for $SO$ and $U$, \emph{Fundamenta Mathematicae} \textbf{174} (2002), 49--78.

\bibitem {realize}\emph{Mathematical Reviews }\#93k:55022, American
Mathematical Society, 1993.

\bibitem {massey}W. Massey, \emph{The mod 2 cohomology of certain Postnikov
systems}, unpublished handwritten manuscript (1978), 20 pages, courtesy of
Haynes Miller.

\bibitem {milnor}J. Milnor, J. Stasheff, \emph{Characteristic Classes}, Annals
of Mathematics Studies 76, Princeton Univ. Press, Princeton, NJ, 1974.

\bibitem {ppw}D. Pengelley, F. Peterson, F. Williams, A global structure
theorem for the mod $2$ Dickson algebras, and unstable cyclic modules over the
Steenrod and Kudo-Araki-May algebras, \emph{Math. Proc. Camb. Phil. Soc.}
\textbf{129 }(2000), 263--275.

\bibitem {pw}D. Pengelley, F. Williams, Sheared algebra maps and operation
bialgebras for mod 2 homology and cohomology, \emph{Transactions of the
American Mathematical Society} \textbf{352} (2000), 1453--1492.

\bibitem {steenrodepstein}N.E. Steenrod, D.B.A. Epstein, \emph{Cohomology
Operations}, Princeton Univ. Press, 1962.

\begin{sloppypar}
\bibitem {stong}R. Stong, Determination of $H^{\ast}(BO(k,\ldots,\infty
),Z_{2})$ and $H^{\ast}(BU(k,\ldots,\infty),Z_{2})$, \emph{Trans. Amer. Math.
Soc.} \textbf{107} (1963), 526--544.
\end{sloppypar}

\bibitem {wilkerson}C. Wilkerson, A primer on the Dickson invariants, in
proc., \emph{Northwestern Homotopy Theory Conference (Evanston, Ill., 1982)},
Contemporary Math. 19 (1983), Amer. Math. Soc., 421--434, as corrected at the
\emph{Hopf Topology Archive},
\texttt{http://hopf.math.purdue.edu/pub/hopf.html}.

\bibitem {wu}W. Wu, Les i-carr\'{e}s dans une vari\'{e}t\'{e}
grassmannienne,\emph{ C.R. Acad. Sci. Paris} \textbf{230} (1950), 918--920.
\end{thebibliography}
\end{document}